\documentclass[journal]{IEEEtran}

\usepackage{amsmath} 
\usepackage{amssymb}  
\usepackage{theorem}
\usepackage{pgf}
\usepackage{mathtools}
\usepackage{tikz}
	\usetikzlibrary{shapes,arrows}

\usepackage{array}
\newcommand{\Ker}{\operatorname{Ker}}
\newcommand{\im}{\operatorname{Im}}

\newcommand{\tr}{\operatorname{trace}}
\newcommand{\diag}{\operatorname{diag}}

\newtheorem{satz}{Theorem}[section]
\sl
\newtheorem{remark}[satz]{Remark}\rm
\newtheorem{lemma}[satz]{Lemma}

\newtheorem{cor}[satz]{Corollary}
\newtheorem{ex}[satz]{Example}

\newcommand{\bprf}{\bf Proof: \nopagebreak  \rm}
\newcommand{\eprf}{\nopagebreak\hspace*{1em} \hfill
$\Box$\noindent\vspace{2ex}\\\phantom{}}

\begin{document}


\title{
Dual pairs of generalized Lyapunov inequalities\\ and balanced
   truncation of stochastic linear systems
}                                

\author{Peter~Benner,
        Tobias~Damm, and Yolanda~Rocio~Rodriguez~Cruz
\thanks{P. Benner is with the Max Planck Institute for Dynamics of Complex Technical Systems, Sandtorstr.\ 1, 39106 Magdeburg, Germany e-mail: benner@mpi-magdeburg.mpg.de.}
\thanks{T. Damm  and Y.R.~Rodriguez Cruz are with University of Kaiserslautern,
Department of Mathematics, 67663 Kaiserslautern, Germany, email:
damm@mathematik.uni-kl.de, rodrigue@mathematik.uni-kl.de}
}


\date{xxxxxx}

\maketitle
\begin{abstract}
 We consider two approaches to balanced truncation of stochastic
 linear systems, which follow from different generalizations of the
 reachability Gram\-ian of deterministic systems. Both preserve mean-square asymptotic
 stability, but only the second leads to a 
 stochastic $H^\infty$-type bound for the approximation error of the truncated
 system. 
\end{abstract}

\section*{Introduction}

Optimization and (feedback) control of dynamical systems is often
computationally infeasible for high dimensional plant models.
Therefore, one tries to reduce the order of the system, so that the
input-output mapping is still computable with sufficient accuracy, but
at considerably smaller cost than for the original system, 
\cite{ObinAnde01, Anto05, BennMehr05, SchiVors08, BaurBenn14}. To
guarantee the desired accuracy, computable error bounds are
required. Moreover, system properties which are relevant in the context of control
system design like asymptotic stability need to be preserved.   It has
long been known that for linear time-invariant (LTI) systems the
method of balanced truncation preserves asymptotic stability and
provides an error bound for the $L^2$-induced input-output norm, that
is the $H^\infty$-norm of the associated transfer function, see \cite{Moor81, PernSilv82}. 
When considering model order reduction of more general system classes,
it is  natural to try to extend this approach. This has been
worked out for descriptor systems in \cite{Styk02}, for
time-varying systems in \cite{ShokSilv83, VerrKail83, SandRant04}, 
for bilinear systems in \cite{AlBaBett94, GrayMesk98, BennDamm11}  and
general nonlinear systems e.g.\ in
\cite{Sche93}. Yet another generaliztion of LTI systems is
obtained considering dynamics driven by noise processes. This leads to
the class of stochastic systems, which have been considered in a system
theoretic context e.g.\ in \cite{Wonh70, HinrPrit98, Damm04}.  
Quite recently, balanced truncation has also
been described for linear stochastic systems of It\^o type in \cite{BennDamm11,
  DammBenn14, RedmBenn14}. 
Already the formulation of the method  leads to two different variants that
are equivalent in the deterministic case, but not so for stochastic systems. 
It is natural to ask which of the above mentioned properties of
balanced truncation also hold for these variants. The aim of this
paper is to answer this question.

Let us first recapitulate balanced truncation for linear
deterministic control systems of the form 
\begin{equation}\label{eq:lti}
 \dot{x} = Ax + Bu,  \quad 
 y       = Cx,\quad \sigma(A)\subset\mathbb{C}_-\;.
\end{equation}
Here $A\in\mathbb{R}^{n\times n}$, $B\in\mathbb{R}^{n\times m}$, $C\in\mathbb{R}^{p\times n}$, and
$x(t)\in\mathbb{R}^n$, $y(t)\in\mathbb{R}^p$ and $u(t)\in\mathbb{R}^m$ are the state,
output, and input of the system, respectively. Moreover $\sigma(A)$ denotes the
spectrum of $A$ and $\mathbb{C}_-$ the open left half complex plane. Let $$\mathcal{L}_A:X\mapsto A^TX+XA$$ denote the
Lyapunov operator and $$\mathcal{L}_A^*:X\mapsto AX+XA^T$$ its adjoint
with respect to the Frobenius inner product. Then
$\sigma(A)\subset\mathbb{C}_-$ if and only if there exists a positive definite solution $X$ of the
 Lyapunov inequality
$\mathcal{L}_A(X)<0$, by Lyapunov's classical stability theorem, see e.g.\ \cite{Gant59b}.

Balanced truncation means truncating a balanced realization. This
realization is obtained by a state space transformation 
%
computed from the Gramians $P$ and $Q$, which solve the dual pair of  \emph{Lyapunov equations}
\begin{subequations}\label{eq:lyap}
  \begin{align}
    \mathcal{L}_A(Q)&=A^T Q + Q A = -C^T C\;,\\ \mathcal{L}_A^*(P)&=A
    P + P A^T = -B B^T\;,
  \end{align}
\end{subequations}
\noindent
or more generally the \emph{inequalities} 
\begin{equation}\label{ineq:lyap}
   \mathcal{L}_A(Q)\le -C^T C\;,\quad
 \mathcal{L}_A^*(P)\le -B B^T\;.  
\end{equation}

These (in)equalities are essential in the
characterization of stability, controllability and observability of 
system \eqref{eq:lti}.
If $\det P\neq 0$, the inequalities \eqref{ineq:lyap}
can be written as
\begin{subequations}\label{eq:lyapequiv}
  \begin{align}
    \mathcal{L}_A(Q)&\le -C^T C\;,\\ \mathcal{L}_A(P^{-1}) &=P^{-1}A +
    A^TP^{-1} \le -P^{-1}B B^TP^{-1}\;.
  \end{align}
\end{subequations}

In the present paper we discuss extensions of \eqref{ineq:lyap} and
\eqref{eq:lyapequiv} for stochastic linear systems. 

As indicated above, the equivalent formulations
\eqref{ineq:lyap} and \eqref{eq:lyapequiv} lead to different
generalizations, if we consider It\^o-type stochastic systems of the form
\begin{equation}\label{eq:ls}
 dx = Ax\,dt+ Nx\,dw + Bu\,dt\;,  \quad 
 y       = Cx\;,
\end{equation}
where $A,B,C$ are as in (\ref{eq:lti}) and
$N\in\mathbb{R}^{n\times n}$. 
System \eqref{eq:ls} is asymptotically mean-square stable (e.g.~\cite{Klei69,
  Khas80, Damm04}), if and only if there
exists a positive definite solution $X$ of the
generalized Lyapunov inequality
\begin{align*}
(\mathcal{L}_A+\Pi_N)(X) = A^TX+XA+N^TXN&<0\;.
\end{align*}
Here $\Pi_N:X\mapsto N^TXN$ and $\Pi_N^*:X\mapsto NXN^T$. This
stability criterion indicates that in the stochastic context the
generalized Lyapunov operator $\mathcal{L}_A+\Pi_N$ takes over the
role of $\mathcal{L}_A$.
Substituting $\mathcal{L}_A$ by $\mathcal{L}_A+\Pi_N$ in  \eqref{ineq:lyap} and \eqref{eq:lyapequiv}, we
obtain two different dual pairs of generalized Lyapunov inequalities.
We call them \emph{type I}:
\begin{subequations}\label{eq:Gram1}
\begin{align}
\label{eq:Gram1.Q}
\displaystyle(\mathcal{L}_A+\Pi_N)(Q)= A^TQ+QA+N^TQN&\le-C^TC\;,\\
\displaystyle (\mathcal{L}_A+\Pi_N)^*(P)=   AP+PA^T+NPN^T&\le-BB^T\;,
\label{eq:Gram1.P}
\end{align}
\end{subequations}
and \emph{type II}:
\begin{subequations}\label{eq:Gram2}
\begin{align}
\label{eq:Gram2.Q}
 (\mathcal{L}_A+\Pi_N)(Q)&=A^TQ+QA+N^TQN\nonumber\\&\le-C^TC\;,\\ 
 (\mathcal{L}_A+\Pi_N)(P^{-1})&= A^TP^{-1}+P^{-1}A+N^TP^{-1}N\nonumber\\&\le-P^{-1}BB^TP^{-1}\;.
 \label{eq:Gram2.P}
\end{align}
\end{subequations}
Note that \eqref{eq:Gram1} corresponds to \eqref{ineq:lyap} in the
sense that $\mathcal{L}_A^*(P)$ has been replaced by
$(\mathcal{L}_A+\Pi_N)^*(P)$, while 
\eqref{eq:Gram2} corresponds to \eqref{eq:lyapequiv}, where
$\mathcal{L}_A(P^{-1})$ has been replaced by $(\mathcal{L}_A+\Pi_N)(P^{-1})$. In general (if $N$ and $P$ do not
commute), the inequalities \eqref{eq:Gram1.P} and
\eqref{eq:Gram2.P} are not equivalent. At first glance it is not
clear which generalization is more appropriate.

If the system is asymptotically mean-square stable and certain
observability and reachability conditions are fulfilled, then 
for both types there are solutions $Q,P>0$. By a suitable state
space-transformation, it is possible to balance the system such that
$Q=P=\Sigma>0$ is diagonal. Consequently, the usual procedure of
balanced truncation can be applied to reduce the order of
\eqref{eq:ls}. For simplicity, let us refer to this as \emph{type I} or \emph{type
II balanced truncation}.

Under natural
assumptions, this reduction preserves mean-square
asymptotic stability.
For type I, this nontrivial fact has been proven in
\cite{BennDamm14c}. Moreover, in \cite{RedmBenn14}, an $H^2$-error
bound has been provided. However, different from the deterministic case,
there is no $H^\infty$-type error bound in terms of the truncated
entries in $\Sigma$. This will be shown in Example \ref{ex:noerrbound}.

In contrast, for type II,   an
$H^\infty$-type error bound has been obtained in \cite{DammBenn14}. 
In the present paper,
as one of our main contributions, we show in Theorem
\ref{satz:Trunc_Stab2} that type II balanced truncation also preserves
mean-square asymptotic stability. The
proof differs significantly  from the one given for type I. Using this result, we are able
to give a more compact proof of the error bound, Theorem
    \ref{thm:NewGramianBound}, which exploits the
stochastic bounded real lemma \cite{HinrPrit98}. 

We illustrate our results by analytical and numerical examples in
Section IV.

\section{Type I balanced truncation}
Consider a stochastic linear control system of It\^o-type
\begin{align}\label{eq:stoch_cont}
      dx&=Ax\,dt+\sum_{j=1}^kN_j x\,dw_j+Bu\,dt\;,\quad 
y= Cx\;,
\end{align}
where $w_j=(w_j(t))_{t\in\mathbb{R}_+}$ are uncorrelated zero mean real 
Wiener processes on a 
probability space $(\Omega ,{\cal F} ,\mu )$ with respect to an
increasing  family 
$({\cal F}_t)_{t\in \mathbb{R}_+}$ of $\sigma$-algebras ${\cal F}_t
\subset{\cal F} $ (e.g.~\cite{Arno74, Oeks98}).\\ To simplify the
notation, we only consider the case $k=1$ and set $w=w_1$, $N=N_1$.
But all results can immediately be generalized for $k>1$. \\
Let $L^2_w(\mathbb{R}_+,\mathbb{R}^q)$ denote the corresponding space of
non-anti\-cipating  stochastic processes $v$ with values in $\mathbb{R}^q$
and norm
\[
\|v(\cdot)\|^2_{L^2_w}:={\cal E}\left(\int_0^\infty\|v(t)\|^2dt\right)<\infty,
\]
where ${\cal E}$ denotes expectation.

Let the homogeneous equation $dx=Ax\,dt+N x\,dw$
be asymptotically mean-square-stable,
i.e.~$\mathcal{E}(\|x(t)\|^2)\stackrel{t\to\infty}\longrightarrow0$,
for all solutions $x$. 

Then, by Theorem \ref{thm:Hans} the equations 
\begin{align*}
\displaystyle A^TQ+QA+N^TQN&=-C^TC\;, \\
 \displaystyle AP+PA^T+NPN^T&=-BB^T\;,
\end{align*}
have unique solutions $Q\ge 0$ and $P\ge 0$. Under suitable
observability and controllability conditions, $Q$ and $P$ are nonsingular.

A similarity
transformation $$(A,N,B,C)\mapsto
(S^{-1}AS,S^{-1}NS,S^{-1}B,CS)$$ of the system implies the 
contragredient transformation  as $$ (Q,P)\mapsto(
S^TQS,S^{-1}PS^{-T})\;.$$ 
Choosing e.g.\ $S=LV\Sigma^{-1/2}$, with Cholesky factoriza\-tions $LL^T=P$, $R^TR=Q$
 and  a singular value
decomposition $RL=U\Sigma V^T$, we obtain $S^{-1}=\Sigma^{-1/2}U^TR$ and
\begin{align*}
  S^{T} Q S &= S^{-1} P S^{-T} =\Sigma =
  \mathop{\mathrm{diag}}(\sigma_1,\ldots,\sigma_n)\;.
\end{align*}
 After suitable partitioning
\begin{align*}
  \Sigma&=\left[
    \begin{array}{cc}
      \Sigma_1&0\\0&\Sigma_2
    \end{array}
  \right],\; S=\left[
    \begin{array}{cc}
      S_1&S_2
    \end{array}
  \right],\; S^{-1}=\left[
    \begin{array}{c}
      T_1\\T_2
    \end{array}
  \right]
\end{align*}
 a truncated system is given in the form $$(A_{11},N_{11},B_1,C_1)=
(T_1AS_1,T_1NS_1,T_1B,CS_1)\;.$$
The following result has been proven in \cite{BennDamm14c}. 
\begin{satz}\label{satz:Trunc_Stab1}
 Let $A,N\in\mathbb{R}^{n\times n}$ satisfy 
 \begin{align*}
   \sigma(I\otimes A+A\otimes I+N\otimes N)\subset\mathbb{C}_-\,.
 \end{align*}
For a
  block-diagonal matrix $\Sigma=\diag(\Sigma_1,\Sigma_2)>0$ with
  $\sigma(\Sigma_1)\cap\sigma(\Sigma_2)=\emptyset$, assume that
  \begin{align*}
    A^T\Sigma+\Sigma A+N^T\Sigma N\le 0\text{ and }
    A\Sigma+\Sigma A^T+N\Sigma N^T\le 0.
  \end{align*}
  Then, with the usual partitioning of $A$ and $N$, we have
  \begin{align*}
    \sigma(I\otimes A_{11}+A_{11}\otimes I+N_{11}\otimes
    N_{11})\subset\mathbb{C}_-\;.
  \end{align*}
\end{satz}
Its implication for mean-square stability of the truncated system is immediate. 
\begin{cor}
  Consider an asymptotically mean square stable stochastic linear system
  \begin{align*}
    dx&=Ax\,dt+Nx\,dw\;.
  \end{align*}
  Assume that a matrix $\Sigma=\diag(\Sigma_1,\Sigma_2)$ is given as
  in Theorem \ref{satz:Trunc_Stab1} and $A$ and $N$ are partitioned
  accordingly.\\
Then the truncated system
\begin{align*}
  dx_r&=A_{11}x_r\,dt+N_{11}x_r\,dw
\end{align*}
is also asymptotically mean square stable.
\end{cor}

If the diagonal entries of $\Sigma_2$
are small, it is expected
that the truncation error is small. In fact this is supported by an $H^2$-error bound
obtained in \cite{RedmBenn14}. Additionally, however, from the
deterministic situation (see \cite{Moor81, Anto05}), one would also hope
for an $H^\infty$-type error bound of the form
\begin{align}\label{eq:alphabound}
   \|y-y_r\|_{L^2_w(\mathbb{R}_+,\mathbb{R}^p)}\stackrel{?}\le\alpha(\tr\Sigma_2)\|u\|_{L^2_w(\mathbb{R}_+,\mathbb{R}^m)}
\end{align}
with some number $\alpha>0$. The following example
shows that no such general $\alpha$ exists.
\begin{ex}\label{ex:noerrbound}\rm
  Let $A=-\left[
  \begin{array}{cc}
    1&0\\0&a^2
  \end{array}
\right]$ with $a>1$, 
$N=\left[
  \begin{array}{cc}
    0&0\\1&0
  \end{array}
\right]$, $B=\left[
  \begin{array}{c}
    1\\0
  \end{array}
\right]$, $C=\left[
  \begin{array}{cc}
    0&1
  \end{array}
\right]$.\\
Solving \eqref{eq:Gram1} with equality, we get   $P=\left[
  \begin{array}{cc}
    \frac12&0\\0&\frac1{4a^2}
  \end{array}
\right]$, $Q=\left[
  \begin{array}{cc}
    \frac1{4a^2}&0\\0&\frac1{2a^2}
  \end{array}
\right]$ with $\sigma(PQ)=\{\frac1{8a^2},\frac1{8a^4}\}$ so that
$\Sigma=\diag(\sigma_1,\sigma_2)$, where $\sigma_1=\frac1{\sqrt8 a}$ and $\sigma_2=\frac1{\sqrt8 a^2}$. The system is
balanced by the transformation $S=\left[
  \begin{array}{cc}
    2a^2&0\\0&1/2
  \end{array}
\right]^{1/4}$.\\ Then $CS=\frac1{2^{1/4}}\left[
  \begin{array}{cc}
    0&1
  \end{array}
\right]$ so that $C_r=0$ for the truncated system of order $1$. 
Thus the output of the reduced system is $y_r\equiv0$, and the
truncation error $\|\mathbb{L}-\mathbb{L}_r\|$ is equal to the
stochastic $H^\infty$-norm (see \cite{HinrPrit98}) of the original system,
\begin{align*}
  \|\mathbb{L}\|=\sup_{x(0)=0, \|u\|_{L^2_w}=1}\|y\|_{L^2_w}\;.
\end{align*} We show now that this norm is
  equal to $\frac1{\sqrt2a}=2a\sigma_2$. Thus, depending on $a$, the ratio of the
  truncation error and $\tr\Sigma_2=\sigma_2$ can be arbitrarily large.

According to the stochastic bounded real lemma, Theorem \ref{thm:sbrl}, 
$\|\mathbb{L}\|$ is the infimum
over all $\gamma$ so that the Riccati inequality
\begin{align}\label{eq:sbrl}
0&< A^TX+XA+N^TXN-C^TC-\frac1{\gamma^2}XBB^TX
\\&=\left[
  \begin{array}{cc}
    -2x_1+x_3-\frac1{\gamma^2}x_1^2&-(a^2+1)x_2-\frac1{\gamma^2}x_1x_2\\
-(a^2+1)x_2-\frac1{\gamma^2}x_1x_2&-2a^2x_3-\frac1{\gamma^2}x_2^2-1
  \end{array}
\right]\nonumber
\end{align}
possesses a solution $X=\left[
  \begin{array}{cc}
    x_1&x_2\\x_2&x_3
  \end{array}
\right]<0$. 

If a given matrix $X$ satisfies this condition, then so
does the same matrix with $x_2$ replaced by $0$. Hence we can assume
that $x_2=0$, and end up with the two conditions
$x_3<-\frac1{2a^2}$ and (after multiplying the upper left entry with $-\gamma^2$)
\begin{align*}
 0&> x_1^2+2\gamma^2x_1-\gamma^2 x_3=(x_1+\gamma^2)^2-\gamma^2(\gamma^2+x_3)\\&>(x_1+\gamma^2)^2-\gamma^2(\gamma^2-\tfrac1{2a^2})\;.
\end{align*}
Thus necessarily $\gamma^2>\frac1{2a^2}$, i.e.\ $\gamma>\frac1{\sqrt2a}$. This already
proves that $\|\mathbb{L}\|\ge \frac1{\sqrt2a}=2a\sigma_2$, which
suffices to disprove the existence of a general bound $\alpha$ in \eqref{eq:alphabound}.
Taking infima, it is easy to show that indeed $\|\mathbb{L}\|=\frac1{\sqrt2a}$.
\end{ex}

\section{Type II balanced truncation}

We now consider the inequalities \eqref{eq:Gram2}.
\begin{lemma}\label{lemma:Gramexists}
  Assume that  $dx=Ax\,dt+N x\,dw$
is asymptotically mean-square-stable. Then  inequality
\eqref{eq:Gram2.P} is solvable with $P>0$.
\end{lemma}
\bprf
 By Theorem \ref{thm:Hans}, for a given $Y<0$, there exists a  $\tilde
P>0$, so that $A^T\tilde P^{-1}+\tilde P^{-1}A+N^T\tilde P^{-1}N=Y$.
Then $P=\varepsilon^{-1}\tilde P$,  for
sufficiently small $\varepsilon>0$, satisfies
\begin{align*}
  A^TP^{-1}+P^{-1}A+N^TP^{-1}N=\varepsilon Y< -\varepsilon^2 \tilde P^{-1}BB^T\tilde P^{-1}
\end{align*}
 so that \eqref{eq:Gram2.P} holds even in the strict form. 
\eprf
It is easy to see that like in the previous section a state space transformation 
\begin{align*}
  (A,N,B,C)\mapsto (S^{-1}AS,S^{-1}NS,S^{-1}B,CS)
\end{align*}
leads to a
contragredient transformation  $Q\mapsto S^TQS$, $ P\mapsto
S^{-1}P S^{-T}$ of the solutions.
 That is, $Q$ and $P$ satisfy \eqref{eq:Gram2.Q} and \eqref{eq:Gram2.P}, if
and only if  $S^TQ S$ and $S^{-1}P S^{-T}$ do so for the transformed data.
As before, we can assume the system to be
balanced with
\begin{align}
  Q&=P=\Sigma=\diag(\sigma_1 I,\ldots,
  \sigma_\nu I)=\left[
    \begin{array}{cc}
      \Sigma_1&\\&\Sigma_2
    \end{array}
\right]\;,\label{eq:defSigma}
\end{align}
where $\sigma_1>\ldots>\sigma_\nu>0$ and
$\sigma(\Sigma_1)=\{\sigma_1,\ldots,\sigma_r\}$, 
$\sigma(\Sigma_2)=\{\sigma_{r+1},\ldots,\sigma_\nu\}$.
Hence,  we will now assume 
(after balancing) that a diagonal matrix
$\Sigma$ as in \eqref{eq:defSigma} is given which satisfies 
\begin{subequations}
  \begin{align}
    A^T\Sigma+\Sigma A+N^T\Sigma N&\le -C^TC\;,\\
    A^T\Sigma^{-1}+\Sigma^{-1}A+ N^T\Sigma^{-1}N&\le
    -\Sigma^{-1}BB^T\Sigma^{-1}\;.
  \end{align}\label{eq:PQ_ineq}
\end{subequations}

Partitioning $A$, $N$, $B$, $C$ like $\Sigma$, we write the system as
\begin{align*}
  dx_1&=(A_{11}x_1+A_{12}x_2)\,dt+(N_{11}x_1+N_{12}x_2)\,dw+B_1u\,dt\\
dx_2&=(A_{21}x_1+A_{22}x_2)\,dt+(N_{21}x_1+N_{22}x_2)\,dw+B_2u\,dt\\
y&=C_1x_1+C_2x_2\;.
\end{align*}
The reduced system obtained by truncation is
\begin{align*}
  dx_r&=A_{11}x_r+N_{11}x_r\,dw+B_1u\,dt\;,\quad y_r=C_1x_r\;.
\end{align*}
The index $r$ is the number of different singular values $\sigma_j$ that
have been kept in the reduced system.
In the following subsections,  we consider matrices 
\begin{align*}
  A=\left[
    \begin{array}{cc}
      A_{11}&A_{12}\\A_{21}&A_{22}
    \end{array}
  \right],\quad N=\left[
    \begin{array}{cc}
      N_{11}&N_{12}\\N_{21}&N_{22}
    \end{array}
  \right]\;, 
\end{align*}
$\Sigma=\diag(\Sigma_1,\Sigma_2)$ as in \eqref{eq:defSigma}, and equations of the form
\begin{subequations}  
\label{eq:CBtilde}
\begin{align}
  A^T\Sigma+\Sigma A+N^T\Sigma N&=-\tilde C^T\tilde C\\
A^T\Sigma^{-1}+\Sigma^{-1}
  A+N^T\Sigma^{-1} N&=-\tilde B\tilde B^T \label{eq:Btilde}
\end{align}
\end{subequations}
with arbitrary right-hand sides $-\tilde
C^T\tilde C\le0$ and $-\tilde B\tilde B^T\le0$.

For convenience, we write out the blocks of these equations explicitly: 
\begin{align}\nonumber
    A_{11}^T\Sigma_1+\Sigma_1A_{11}+N_{11}^T&\Sigma_1N_{11}
\\&=-N_{21}^T\Sigma_2
  N_{21}-\tilde C_1^T\tilde C_1\label{eq:allcomps1}
 \\[2mm]\nonumber
A_{12}^T\Sigma_1+\Sigma_2A_{21}+N_{12}^T&\Sigma_1N_{11}
\\&=-N_{22}^T\Sigma_2
  N_{21}-\tilde C_2^T\tilde C_1\label{eq:allcomps2}\\[2mm]
\nonumber
A_{22}^T\Sigma_2+\Sigma_2A_{22}+N_{22}^T&\Sigma_2N_{22}
\\&=-N_{12}^T\Sigma_1
  N_{12}-\tilde C_2^T\tilde C_2\label{eq:allcomps3}\\[2mm]
\nonumber
A_{11}^T\Sigma^{-1}_1+\Sigma^{-1}_1A_{11}+N_{11}^T&\Sigma^{-1}_1N_{11}
\\&=-N_{21}^T\Sigma^{-1}_2
  N_{21}-\tilde B_1\tilde B_1^T\label{eq:allcomps4}
\\[2mm]\nonumber
A_{12}^T\Sigma^{-1}_1+\Sigma^{-1}_2A_{21}+N_{12}^T&\Sigma^{-1}_1N_{11}
\\&=-N_{22}^T\Sigma^{-1}_2
  N_{21}-\tilde B_2\tilde B_1^T\label{eq:allcomps5}\\[2mm]
\nonumber
A_{22}^T\Sigma^{-1}_2+\Sigma^{-1}_2A_{22}+N_{22}^T&\Sigma^{-1}_2N_{22}
\\&=-N_{12}^T\Sigma^{-1}_1
  N_{12}-\tilde B_2\tilde B_2^T\label{eq:allcomps6}
\end{align}

\subsection{Preservation of asymptotic stability}

The following theorem is the main new result of this paper. 
\begin{satz}\label{satz:Trunc_Stab2}
  Let $A$ and $N$ be given such that
  \begin{align}
    \sigma(I\otimes A+A\otimes I+N\otimes N)\subset\mathbb{C}_-\;.\label{eq:Assumedstability}   
  \end{align}
Assume further that for a
  block-diagonal matrix $\Sigma=\diag(\Sigma_1,\Sigma_2)>0$ with
  $\sigma(\Sigma_1)\cap\sigma(\Sigma_2)=\emptyset$, we have
  \begin{subequations}\label{eq:Sigma_ineq}
    \begin{align}
      A^T\Sigma+\Sigma A+N^T\Sigma N&\le 0\quad\text{ and }\\
      A^T\Sigma^{-1}+\Sigma^{-1} A+N^T\Sigma^{-1} N&\le 0\;.
    \end{align}
  \end{subequations}
  Then, with the usual partitioning of $A$ and $N$, we have
  \begin{align}\label{eq:Trunc_Stab2}
    \sigma(I\otimes A_{11}+A_{11}\otimes I+N_{11}\otimes
    N_{11})\subset\mathbb{C}_-\;.
  \end{align}
\end{satz}
Again we have an immediate interpretation in terms of mean-square stability of the truncated system. 
\begin{cor}
  Consider an asymptotically mean square stable stochastic linear system
  \begin{align*}
    dx&=Ax\,dt+Nx\,dw\;.
  \end{align*}
  Assume that a matrix $\Sigma=\diag(\Sigma_1,\Sigma_2)$ is given as
  in Theorem \ref{satz:Trunc_Stab2} and $A$ and $N$ are partitioned
  accordingly.\\
Then the truncated system
\begin{align*}
  dx_r&=A_{11}x_r\,dt+N_{11}x_r\,dw
\end{align*}
is also asymptotically mean square stable.
\end{cor}

{\bf Proof of Theorem \ref{satz:Trunc_Stab2}: \nopagebreak  \rm} 
Note that the inequalities
\eqref{eq:Sigma_ineq}
are equivalent to the equations 
\eqref{eq:allcomps1} -- \eqref{eq:allcomps6} with appropriate
right-hand sides $-\tilde C^T\tilde C$ and $-\tilde
B\tilde B^T$.\\
By way of contradiction, we assume that \eqref{eq:Trunc_Stab2} does not hold.
Then by Theorem \ref{thm:KreinRutman}, there exist $V\ge 0$, $V\neq 0$, $\alpha\ge 0$ such that
\begin{align}\label{eq:evecV}
  A_{11}V+VA_{11}^T+N_{11}VN_{11}^T&=\alpha V\;.
\end{align}
Taking the scalar product of the equation \eqref{eq:allcomps1} with $V$, we obtain $0\ge
\alpha\tr(\Sigma_1V)$ whence $\alpha=0$ and  $\tilde C_1V=0$,
$N_{21}V=0$  by
Corollary \ref{cor:semidefY}. Hence
\begin{align}
  \label{eq:scprodV}
  \left(A_{11}^T\Sigma_1+\Sigma_1A_{11}+N_{11}^T\Sigma_1N_{11} \right) V&=0\;.
\end{align}
Analogously, we have $\tilde
B_1^TV=0$ by \eqref{eq:allcomps2}.\\
In particular, from $N_{21}V=0$, we get
\begin{align*}
&(\mathcal{L}_A^*+\Pi_N^*)\left(\left[\begin{array}{cc}V&0\\0&0\end{array}\right]\right)\\
&=\left[
  \begin{array}{cc}
    A_{11}V+VA_{11}^T+N_{11}VN_{11}^T&  VA_{21}^T+N_{11}VN_{21}^T \\
A_{21}V+N_{21}VN_{11}^T&N_{21}VN_{21}^T
  \end{array}
\right]\\&=\left[
  \begin{array}{cc}
    0&  VA_{21}^T \\
A_{21}V&0
  \end{array}
\right]\;.
\end{align*}
We will show that $A_{21}V=0$, which implies
\begin{align}
  0\in\sigma(I\otimes A+A\otimes I+N\otimes N)\label{eq:contradict}
\end{align}
in contradiction to \eqref{eq:Assumedstability},
  and thus finishes the proof.

We first show that $\im V$ is invariant under $A_{11}$ and $N_{11}$.
To this end let $Vz=0$. Then by \eqref{eq:evecV},
  \begin{align*}
    0&=z^T\left(
      A_{11}V+VA_{11}^T+N_{11}VN_{11}^T\right)z=z^TN_{11}VN_{11}^Tz\;,
  \end{align*}
  whence also $VN_{11}^Tz=0$, i.e.\ $N_{11}^Tz\in\Ker V$. From this,
  we have
  \begin{align*}
    0&=\left(
      A_{11}V+VA_{11}^T+N_{11}VN_{11}^T\right)z=VA_{11}^Tz\;,
  \end{align*}
  implying $A_{11}^Tz\in\Ker V$. Thus $A_{11}^T\Ker V\subset \Ker V$
  and $N_{11}^T\Ker V\subset \Ker V$.

  Since $\Ker V=(\im V)^\bot$, it follows further that $\im V$
  is invariant under $A_{11}$ and $N_{11}$.

Let $V=V_1V_1^T$, where $V_1$ has full column rank, i.e.\ $\det
V_1^TV_1\neq 0$. Then by the
invariance, there exist square matrices $X$ and $Y$, such that
\begin{align*}
  A_{11}V_1=V_1X\quad\text{ and }\quad N_{11}V_1=V_1Y\;.
\end{align*}
It follows that
\begin{align*}
  0&=A_{11}V_1V_1^T+V_1V_1^TA_{11}^T+N_{11}V_1V_1^TN_{11}^T\\&=V_1(X+X^T+YY^T)V_1^T\;,
\end{align*}
whence $X+X^T+YY^T=0$. Moreover, from \eqref{eq:scprodV},
we get
\begin{align}\nonumber
  A_{11}^T\Sigma_1V_1&=-\Sigma_1A_{11}V_1-N_{11}^T\Sigma_1N_{11}V_1\\&=-\Sigma_1V_1X-N_{11}^T\Sigma_1V_1Y\;.
\label{eq:A11SigmaV1}
\end{align}
Using this substitution in the following computation, we obtain
\begin{align}\nonumber
  0&\ge V_1^T\Sigma_1^2\left(A_{11}^T\Sigma_1^{-1}+\Sigma_1^{-1}A_{11}+N_{11}^T\Sigma_1^{-1}N_{11}\right)\Sigma_1^2V_1\\\nonumber
&=
V_1^T\Sigma_1^2(A_{11}^T\Sigma_1V_1)+(A_{11}^T\Sigma_1V_1)^T\Sigma_1^2V_1\\\nonumber&\phantom{==}+V_1^T\Sigma_1^2N_{11}^T\Sigma_1^{-1}N_{11}\Sigma_1^2V_1\\\label{eq:Sigma1hoch3}
&=-V_1^T\Sigma_1^3V_1X-X^TV_1^T\Sigma_1^3V_1\\\nonumber
&\phantom{==}-V_1^T\Sigma_1^2N_{11}^T\Sigma_1V_1Y-Y^TV_1^T\Sigma_1N_{11}\Sigma_1^2V_1\\\nonumber
&\phantom{==}+V_1^T\Sigma_1^2N_{11}^T\Sigma_1^{-1}N_{11}\Sigma_1^2V_1\;.
\end{align}
We will show that the right hand side has nonnegative trace. This then
implies that the whole term vanishes. Note that
\begin{align*}
\tr (Y^TV_1^T\Sigma_1^3V_1Y)&=\tr (V_1^T\Sigma_1^3V_1YY^T)\\&=\tr (-V_1^T\Sigma_1^3V_1(X+X^T))\\
  &=\tr (-V_1^T\Sigma_1^3V_1X-X^TV_1^T\Sigma_1^3V_1)\;.
\end{align*}
Taking the trace in \eqref{eq:Sigma1hoch3}, we have
\begin{align*}
  0&\ge\tr
 \Big(Y^TV_1^T\Sigma_1^3V_1Y-V_1^T\Sigma_1^2N_{11}^T\Sigma_1V_1Y
\\&\phantom{xxxxxx}-Y^TV_1^T\Sigma_1N_{11}\Sigma_1^2V_1
+V_1^T\Sigma_1^2N_{11}^T\Sigma_1^{-1}N_{11}\Sigma_1^2V_1\Big)\\
&=\tr\left[
  \begin{array}{c}
    V_1Y\\V_1
  \end{array}
\right]^T M
\left[
  \begin{array}{c}
    V_1Y\\V_1
  \end{array}
\right] \;.
\end{align*}
where
\begin{align*}
M&=     \left[
  \begin{array}{cc} \Sigma_1^3&-\Sigma_1N_{11}\Sigma_1^{2}\\-\Sigma_1^2N_{11}^T\Sigma_1&\Sigma_1^2N_{11}^T\Sigma_1^{-1}N_{11}\Sigma_1^2
  \end{array}
\right]\;.
\end{align*}
The matrix $M$ is positive semidefinite,
because the upper left block is positive definite, and the
corresponding Schur complement
\begin{align*} \Sigma_1^2N_{11}^T\Sigma_1^{-1}N_{11}\Sigma_1^2-\Sigma_1^2N_{11}^T\Sigma_1\Sigma_1^{-3}\Sigma_1N_{11}\Sigma_1^{2}&=0
\end{align*}
vanishes. Hence 
\begin{align*}
  \left[
    \begin{array}{cc}
      \Sigma_1^3&-\Sigma_1N_{11}\Sigma_1^{2}\\
-\Sigma_1^2N_{11}^T\Sigma_1&\Sigma_1^2N_{11}^T\Sigma_1^{-1}N_{11}\Sigma_1^2
    \end{array}
  \right]\left[
    \begin{array}{c}
      V_1Y\\V_1
    \end{array}
  \right]&=0
\end{align*}
implying via the first block row that $N_{11}\Sigma_1^2 V_1=\Sigma_1^2 V_1Y$. From
\eqref{eq:Sigma1hoch3}, using also \eqref{eq:A11SigmaV1} again, we
thus have
\begin{align*}
0&=
\left(A_{11}^T\Sigma_1^{-1}+\Sigma_1^{-1}A_{11}+N_{11}^T\Sigma_1^{-1}N_{11}\right)\Sigma_1^2V_1\\
&=-\Sigma_1V_1X-N_{11}^T\Sigma_1V_1Y+\Sigma_1^{-1}A_{11}\Sigma_1^2V_1+N_{11}^T\Sigma_1
V_1Y\\
&=-\Sigma_1V_1X+\Sigma_1^{-1}A_{11}\Sigma_1^2V_1\;,
\end{align*}
i.e.\  $A_{11}\Sigma_1^2V_1=\Sigma_1^2V_1X$. It follows that for
arbitrary $k\in\mathbb{N}$, the
eigenvector $V$ in \eqref{eq:evecV} can be replaced by
$$\Sigma_1^{2k}V\Sigma_1^{2k}=\Sigma_1^{2k}V_1V_1^T\Sigma_1^{2k}$$ because
\begin{align*}
0&=\Sigma_1^2V_1\left(X+X^T+YY^T\right)V_1^T\Sigma_1^2\\
&=  A_{11}\left(\Sigma_1^2V_1V_1^T\Sigma_1^2\right)+\left(\Sigma_1^2V_1V_1^T\Sigma_1^2\right)A_{11}^T
\\&\phantom{==}+N_{11}\left(\Sigma_1^2V_1V_1^T\Sigma_1^2\right)N_{11}^T\;.
\end{align*}
Induction leads to
\begin{align*}
0&=
A_{11}\left(\Sigma_1^{2k}V_1V_1^T\Sigma_1^{2k}\right)+\left(\Sigma_1^{2k}V_1V_1^T\Sigma_1^{2k}\right)A_{11}^T\\
&\phantom{==}+ N_{11}\left(\Sigma_1^{2k}V_1V_1^T\Sigma_1^{2k}\right)N_{11}^T\;.
\end{align*}
As above, we conclude that $N_{21}\Sigma_1^{2k}V_1=0$,
$\tilde C_1\Sigma_1^{2k}V_1=0$, and $\tilde B_1^T\Sigma_1^{2k}V_1=0$. 
Multiplying \eqref{eq:allcomps2} with $\Sigma_1^{2(k-1)}V_1$
and \eqref{eq:allcomps5} with $\Sigma_1^{2k}V_1$, we get
\begin{align*}
 A_{12}^T\Sigma_1^{2k-1}V_1+\Sigma_2A_{21}\Sigma_1^{2(k-1)}V_1+N_{12}^T\Sigma_1^{2k-1}V_1Y&=0\;,\\
A_{12}^T\Sigma_1^{2k-1}V_1+\Sigma_2^{-1}A_{21}\Sigma_1^{2k}V_1+N_{12}^T\Sigma_1^{2k-1}V_1Y&=0\;.
\end{align*}
Hence (after multiplication with $\Sigma_2$), for all $k\ge 1$, we have
\begin{align*}
  \Sigma_2^2A_{21}\Sigma_1^{2(k-1)}V_1&=-\Sigma_2\left(
    A_{12}^T\Sigma_1^{2k-1}V_1+N_{12}^T\Sigma_1^{2k-1}V_1Y \right)\\&=A_{21}\Sigma_1^{2k}V_1\;.
\end{align*}
Applying this identity repeatedly, we get
$$A_{21}\Sigma_1^{2k}V_1=\Sigma_2^{2k}A_{21}V_1 \quad\text{ for all $k\in\mathbb{N}$. }$$
If $\mu$ is the minimal polynomial of $\Sigma_1^2$, then
$\sigma(\Sigma_1)\cap\sigma(\Sigma_2)=\emptyset$ implies
$\det\mu(\Sigma_2^2)\neq0$ and
\begin{align*}
  0&=A_{21}\mu(\Sigma_1^2)V_1=\mu(\Sigma_2^2)A_{21}V_1\;,
\end{align*}
whence $A_{21}V_1=0$ and also $A_{21}V=0$. Hence we obtain the
contradiction \eqref{eq:contradict}.
\eprf

\subsection{Error estimate}
The following theorem has been proven in \cite{DammBenn14} using
LMI-techniques. Exploiting the stability result in the previous
subsection, we can give a slightly more compact proof based on the
stochastic bounded real lemma, Theorem \ref{thm:sbrl0}. 

\begin{satz}\label{thm:NewGramianBound} 
 Let $A$ and $N$ satisfy
 \begin{align*}
   \sigma(I\otimes A+A\otimes I+N\otimes N)\subset\mathbb{C}_-\;.
 \end{align*}
Assume furthermore that for 
  $\Sigma=\diag(\Sigma_1,\Sigma_2)>0$ with
  $\Sigma_2=\diag(\sigma_{r+1} I,\ldots,\sigma_\nu I)$
  and  $\sigma(\Sigma_1)\cap\sigma(\Sigma_2)=\emptyset$,  the following
  Lyapunov inequalities hold,
  \begin{align*}
    A^T\Sigma+\Sigma A+N^T\Sigma N&\le -C^TC\;,\\
    A^T\Sigma^{-1}+\Sigma^{-1} A+N^T\Sigma^{-1} N&\le -\Sigma^{-1}BB^T\Sigma^{-1}\;.
  \end{align*}
If $x(0)=0$ and $x_r(0)=0$, then for all $T>0$, it holds that
 $$\displaystyle \|y-y_r\|_{L^2_w([0,T])}\le2
 (\sigma_{r+1}+\ldots+\sigma_\nu)\|u\|_{L^2_w([0,T])}\;.$$
\end{satz}
\bprf We adapt a proof for deterministic systems e.g.\
\cite[Theorem 7.9]{Anto05}. In the central argument we treat the case
where  
$\Sigma_2=\sigma_\nu I$ and show that
\begin{align}\label{eq:allpass_bound}
  \|y-y_{\nu-1}\|_{L^2_w[0,T]}\le 2\sigma_\nu \|u\|_{L^2_w[0,T]}\;.
\end{align}
From \eqref{eq:allcomps1} and \eqref{eq:allcomps4}, we can see that also 
\begin{align*}
    A_{11}^T\Sigma_1+\Sigma_1 A_{11}+N_{11}^T\Sigma_1 N_{11}&\le -C_1^TC_1\;,\\
    A_{11}^T\Sigma_1^{-1}+\Sigma_1^{-1} A_{11}+N_{11}^T\Sigma_1^{-1} N_{11}&\le -\Sigma_1^{-1}B_1B_1^T\Sigma_1^{-1}\;.
  \end{align*}
Hence we can repeat the above argument 
to remove
$\sigma_{\nu-1},\ldots,\sigma_{r+1}$ successively. By the triangle inequality we
find that 
\begin{align*}
  \|y-y_r\|_{L^2_w[0,T]}&\le \sum_{j=r}^{\nu-1}\|y_{j+1}-y_j\|_{L^2_w[0,T]}\\
&  \le 2(\sigma_\nu+\ldots+\sigma_{r+1}) \|u\|_{L^2_w[0,T]}\;.
\end{align*}
 which then concludes the proof.\\
To prove \eqref{eq:allpass_bound}, we make use of the stochastic
bounded real lemma.
In the following let $r=\nu-1$ and consider the error system
defined by
\begin{align*}
  dx_e&=A_ex_e\,dt+N_ex_e\,dw+B_eu\,dt\;,\\ y_e&=C_ex_e=y-y_r\;,
\end{align*}
where 
\begin{align*}
  x_e&=\left[
    \begin{array}{c}
      x_1\\x_2\\x_r
    \end{array}
\right],\quad A_e=\left[
    \begin{array}{ccc}
      A_{11}&A_{12}&0\\A_{21}&A_{22}&0\\0&0&A_{11}
    \end{array}
\right],\\ N_e&=\left[
    \begin{array}{ccc}
      N_{11}&N_{12}&0\\N_{21}&N_{22}&0\\0&0&N_{11}
    \end{array}
\right],\quad B_e=\left[
    \begin{array}{c}
      B_1\\B_2\\B_1
    \end{array}
\right],\\ C_e&=\left[
    \begin{array}{ccc}
      C_1&C_2&-C_1
    \end{array}
\right]\;.
\end{align*}
Applying the state space transformation
\begin{align*}
  \left[
    \begin{array}{c}
      \tilde x_1\\\tilde x_2\\\tilde x_r
    \end{array}
\right]&=\left[
    \begin{array}{c}
      x_1-x_r\\x_2\\x_1+x_r
    \end{array}
\right]=\underbrace{\left[
    \begin{array}{ccc}
      I_r&0&-I_r\\0&I_{n-r}&0\\I_r&0&I_r
    \end{array}
\right]}_{=S^{-1}}\left[
    \begin{array}{c}
      x_1\\x_2\\x_r
    \end{array}
\right], 
\end{align*}
we obtain the transformed system
\begin{align*}
  \tilde A_e&=S^{-1}A_eS=\left[
    \begin{array}{ccc}
      A_{11}&A_{12}&0\\\tfrac12A_{21}&A_{22}&\tfrac12A_{21}\\0&A_{12}&A_{11}
    \end{array}
\right]\;,\\\tilde N_e&=S^{-1}N_eS=\left[
    \begin{array}{ccc}
      N_{11}&N_{12}&0\\\tfrac12N_{21}&N_{22}&\tfrac12N_{21}\\0&N_{12}&N_{11}
    \end{array}
\right]\;,\\ \tilde
  B_e&=S^{-1}B\left[
    \begin{array}{c}
      0\\B_2\\2B_1
    \end{array}
\right]\;,\quad\tilde C_e=C_eS=\left[
    \begin{array}{ccc}
      C_1&C_2&0
    \end{array}
\right]\;.
\end{align*}
By Theorem \ref{thm:sbrl0}, we have
$\|\mathbb{L}_e\|\le 2\sigma_\nu$, if the Riccati inequality  
\begin{align}\nonumber
\mathcal{R}_\gamma(X)&= \tilde A_e^TX+X\tilde A_e+\tilde N_e^TX\tilde N_e+\tilde
C_e^T\tilde C_e\\&\phantom{==}+\frac1{4\sigma_\nu^2}X\tilde B_e\tilde
B_e^TX\le 0
\label{eq:sbrl_error}
\end{align}
possesses a solution $X\ge 0$. We will show now that the 
block-diagonal matrix
\begin{align*}
 X=\diag(\Sigma_1,2\Sigma_2,\sigma_\nu^2
  \Sigma_1^{-1})=\diag(\Sigma_1,2\sigma_\nu I,\sigma_\nu^2
  \Sigma_1^{-1})>0
\end{align*}
satisfies \eqref{eq:sbrl_error}. Partitioning
  $\mathcal{R}_{\sigma_\nu}(X)=\left[
    \begin{array}{ccc}
      R_{11}&R_{21}^T&R_{31}^T\\R_{21}&R_{22}&R_{32}^T\\R_{31}&R_{32}&R_{33}
    \end{array}
\right]$, 
we have 
\begin{align*}
  R_{11}&=A_{11}^T\Sigma_1+\Sigma_1A_{11}+N_{11}^T\Sigma_1N_{11}+\frac{\sigma_\nu}2N_{21}^TN_{21}+C_1^TC_1\\
&=A_{11}^T\Sigma_1+\Sigma_1A_{11}+N_{11}^T\Sigma_1N_{11}+
N_{21}^T\Sigma_2N_{21}+C_1^TC_1\\&\phantom{==}-\frac{\sigma_\nu}2N_{21}^TN_{21}\\
R_{21}&=A_{12}^T\Sigma_1+{\sigma_\nu}
A_{21}+N_{12}^T\Sigma_1N_{11}+{\sigma_\nu} N_{22}^TN_{21}+C_2^TC_1\\
R_{31}&=\frac{\sigma_\nu}2 N_{21}^TN_{21}\\
R_{22}&=2{\sigma_\nu}(A_{22}^T+A_{22}+N_{22}^TN_{22})+N_{12}^T\Sigma_1N_{12}\\&\phantom{==}+
\sigma_\nu^2N_{12}^T\Sigma_1^{-1}N_{12}+C_2^TC_2+B_2B_2^T\\
&=A_{22}^T\Sigma_2+\Sigma_2A_{22}+N_{22}^T\Sigma_2N_{22}+N_{12}^T\Sigma_1N_{12}+C_2^TC_2\\
&\phantom{==} +\sigma_\nu^2(A_{22}^T\Sigma_2^{-1}+\Sigma_2^{-1}A_{22}+
N_{22}^T\Sigma_2^{-1}N_{22}\\
&\phantom{==}+N_{12}^T\Sigma_1^{-1}N_{12}+\Sigma_2^{-1}B_2B_2^T\Sigma_2^{-1})\\
R_{32}&=\sigma_\nu^2(\Sigma_1^{-1}A_{12}+N_{11}^T\Sigma_1^{-1}N_{12})+{\sigma_\nu}
(A_{21}^T+N_{21}^TN_{22})\\
&\phantom{==}+{\sigma_\nu}\Sigma_1^{-1}B_1B_2^T\\
&=\sigma_\nu^2(\Sigma_1^{-1}A_{12}+N_{11}^T\Sigma_1^{-1}N_{12}+A_{21}^T\Sigma_2^{-1}+N_{21}^T\Sigma_2^{-1}N_{22}\\
&\phantom{==}+\Sigma_1^{-1}B_1B_2^T\Sigma_2^{-1})\\
R_{33}&=\sigma_\nu^2(A_{11}^T\Sigma_1^{-1}+\Sigma_1^{-1}A_{11}+N_{11}^T\Sigma_1^{-1}N_{11})+\frac{\sigma_\nu}2
N_{21}^TN_{21}\\
&\phantom{==}+\sigma_\nu^2\Sigma_1^{-1}B_1B_1^T\Sigma_1^{-1}\\
&=\sigma_\nu^2(A_{11}^T\Sigma_1^{-1}+\Sigma_1^{-1}A_{11}+N_{11}^T\Sigma_1^{-1}N_{11}
\\
&\phantom{==}+
\Sigma_1^{-1}B_1B_1^T\Sigma_1^{-1}+N_{21}^T\Sigma_2^{-1}N_{21})-\frac{\sigma_\nu}2N_{21}^TN_{21}
\end{align*}
With the permutation matrix $J=\left[
  \begin{array}{cc}
    0&I\\I&0
  \end{array}
\right]$ we define 
\begin{align*}
  M&=J(A^T\Sigma^{-1}+\Sigma^{-1} A+N^T\Sigma^{-1} N+\Sigma^{-1}BB^T\Sigma^{-1})J\;,
\end{align*}
where $M\le 0$ by \eqref{eq:Btilde}.
Using \eqref{eq:allcomps1} -- \eqref{eq:allcomps6}, we have
\begin{align*}
 & \mathcal{R}_{\sigma_\nu}(X)=\left[
  \begin{array}{c|c}
    A^T\Sigma+\Sigma A+N^T\Sigma N+C^TC&0\\\hline0&0
  \end{array}
\right]\\
&-\frac{\sigma_\nu}2\left[
    \begin{array}{c}
      N_{21}^T\\0\\-N_{21}^T
    \end{array}
\right]\left[
    \begin{array}{c}
      N_{21}^T\\0\\-N_{21}^T
    \end{array}
\right]^T+\sigma_\nu^2\left[
  \begin{array}{c|c}
0&0\\\hline0&
    M
  \end{array}
\right]\le 0\;,
\end{align*}
which is inequality \eqref{eq:sbrl_error}.
\eprf

\begin{ex}\rm \label{ex:noerror_cont}
  Let the system $(A,N,B,C)$ and $Q$ be as in Example~\ref{ex:noerrbound}. The matrix
\begin{align*}
  P=\left[
    \begin{array}{cc}
1+\sqrt{1-p}&0\\0&p
\end{array}
\right]^{-1}>0\;,\text{ where } 0<p\le 1\;,
\end{align*}
satisfies inequality \eqref{eq:Gram2.P}.
As in Example \ref{ex:noerrbound}, we have $\mathbb{L}_r=0$ for the
corresponding reduced system of
order $1$, so that the truncation error again is $\frac1{\sqrt2a}$,
independently of $p\in\,]0,1]$.

On the other hand we have 
  \begin{align*}
   \sigma_2^2&= \min\sigma(PQ)=\frac1{4a^2(1+\sqrt{1-p})}\le \frac1{8a^2}\;,
  \end{align*}
with equality for $p\to 0$. Theorem \ref{thm:NewGramianBound} thus
gives the sharp error bound $2\sigma_2=\frac1{\sqrt2a}$. Note, that
there is no $P>0$ satisfying the \emph{equation} \eqref{eq:Gram2.P}.
\end{ex}
  The previous example illustrates the problem of optimizing over all solutions of
  inequality \eqref{eq:Gram2.P}.
\section{Numerical Examples}
To compare the reduction methods we need to compute $Q, P$ from
\eqref{eq:Gram1} or \eqref{eq:Gram2}.
Instead of the inequalities \eqref{eq:Gram1.Q},
\eqref{eq:Gram1.P}, \eqref{eq:Gram2.Q} we can consider the corresponding equations, for which quite efficient  algorithms
have been developed recently, e.g.\ \cite{Damm08, BennBrei13,
  KresSirk14, ShanSimo14}. These also allow for a low-rank
approximation of the solutions. In
contrast we cannot replace  \eqref{eq:Gram2.P} by the corresponding
equation, because this may not be solvable (see Example
\ref{ex:noerror_cont}). Even worse, we do not have any solvability or uniqueness
criteria nor reliable algorithms.

 Therefore, in general, we have to work with the inequality
 \eqref{eq:Gram2.P}, which is solvable according to Lemma
 \ref{lemma:Gramexists}, but of course not uniquely solvable.

    In view of our application, we aim at a solution $P$ of \eqref{eq:Gram2.P}, so that
    (some of) the eigenvalues of $PQ$ are particularly small,
    since they provide the error bound. Choosing a matrix $Y<0$ and a very small
    $\varepsilon$ along the lines of the proof of Lemma
 \ref{lemma:Gramexists} can be contrary to this aim.
    Hence some optimization over all solutions of \eqref{eq:Gram2.P}
    is required. 

Note also that a matrix $P>0$ satisfies
    \eqref{eq:Gram2.P}, if and only if it satisfies the linear matrix
    inequality (LMI)
    \begin{align}\label{eq:lmi}
      \left[
        \begin{array}{cc}
          PA^T+AP+BB^T&PN^T\\NP&-P
        \end{array}
      \right]&\le 0\;.
    \end{align}
Thus, LMI optimal solution techniques are applicable. However, their
complexity will be prohibitive for large-scale problems. Therefore
further research for alternative methods to solve \eqref{eq:Gram2.P}
 adequately is required. 

By $\mathbb{L}$ and
$\mathbb{L}_r$, we always denote the original and the $r$-th order
approximated system. The stochastic $H^\infty$-type norm
$\|\mathbb{L}-\mathbb{L}_r\|$ is computed by a binary search of the
infimum of all $\gamma$ such that the Riccati inequality \eqref{eq:sbrl}  is solvable. The latter is solved via a Newton iteration as in \cite{Damm04}.
Finally, the Lyapunov equations
\eqref{eq:lyap} are solved by preconditioned Krylov subspace methods
described in \cite{Damm08}. \\
Unfortunately, for small $\gamma$, i.e.\ for small
approximation errors, this method of computing the error runs into numerical problems,
because  \eqref{eq:sbrl} contains the term $\gamma^{-2}$. This
apparently leads to cancellation phenomena in the Newton iteration, if
e.g.\ $\gamma<10^{-7}$.  Therefore we mainly concentrate on cases where
the error is larger, that is we make $r$ sufficiently small.

\subsection{Type II can be better than type I}
In many examples we observe that type II reduction gives a valid error
bound, but the approximation error still is better with type I. This,
however, is not always true, as the example
\begin{align*}
  (A,N,B,C^T)=\left(\left[
    \begin{array}{rr}
-1&1\\0&-1
\end{array}
\right],\left[
    \begin{array}{rr}
0&0\\1&0
\end{array}
\right],
\left[
    \begin{array}{r}
0\\3
\end{array}
\right],
\left[
    \begin{array}{r}
3\\0
\end{array}
\right]\right)
\end{align*}
shows. It can easily be verified that the type I Lyapunov equations \eqref{eq:Gram1}
are solved by
\begin{align*}
  Q=\left[
    \begin{array}{rr}
6&3\\3&3
\end{array}
\right]\text{ and }
P=\left[
    \begin{array}{rr}
3&3\\3&6
\end{array}
\right]\;.
\end{align*}
The type II inequalities \eqref{eq:Gram2} are e.g.\ solved by 
\begin{align*}
  Q=\left[
    \begin{array}{rr}
6&3\\3&3
\end{array}
\right]\text{ and }
P=\left[
    \begin{array}{rr}
8&0\\0&12
\end{array}
\right]\;.
\end{align*}
Reduction to order $r=1$ gives the following error bounds and
approximation errors for both types:
\begin{center}
  \begin{tabular}{c||c|c}
    &$\sigma_2$&$\|\mathbb{L}-\mathbb{L}_{1}\|$\\\hline
    I &$2.4853$&$3.9647$\\\hline
    II &$6.9282$&$3.5614$
  \end{tabular}
\end{center}
As we see, the type I approximation error is larger than both the
truncated singular value and the type II approximation error.

\subsection{An electrical ladder network with perturbed inductance}
\label{sec:an-electrical-ladder}
As our first example with a physical background, we take up the electrical ladder network described in \cite{GugeAnto04}, consisting of $n/2$ sections with a
capacitor $\tilde C$, inductor $\tilde L$  and two resistors $R$ and $\tilde R$ as depicted
here.

\begin{center}
     \tikzstyle{blockh} = [draw, rectangle, minimum height=2em, minimum
     width=2.5em] 
     \tikzstyle{blockv} = [draw, rectangle, minimum height=3em,
     minimum width=1.5em] 
\begin{tikzpicture}
  \draw[dashed] (0,0) rectangle(4.5,3.8);
  \node at (1,3) [blockh](R){$R=0.1$};
  \node at (3,3) [blockh](L){$\tilde L=0.1$};
  \node at (2,1.5) [blockv](C){{\rotatebox{-90}{$\tilde C=0.1$}}};
  \node at (4,1.5) [blockv](R2){{\rotatebox{-90}{$\tilde R=1$}}};
  \node at (-.8,1.75) {$V$};
  \node at (.6,.8) {$I$};
  \draw[<-] (-.5,2.8) to (-.5,.7);
  \draw[<-] (.5,.5) to (5,.5);
  \draw (-.5,3) to (R) to (L) to (5,3);
  \draw(2,3) to (C) to (2,.5);
  \draw(4,3) to (R2) to (4,.5);
  \draw(-.5,.5) to (.5,.5);
\end{tikzpicture}

\end{center}
 But following e.g.\ \cite{UgriPete99}, we
assume that the inductance $\tilde L$ is subject to stochastic perturbations.
For simplicity, we replace the inverse $\tilde L^{-1}$ formally by $L^{-1}+\dot
w$ in all sections. Here $L=0.1$ and $\dot w$ is white noise of a
certain intensity $\sigma$, where we set $\sigma=1$. E.g.\ for $n=6$,
we have the system matrices
{\small
\begin{align*}
  A&=\left[\begin{array}{cccccc}
 \frac{-1}{\tilde C R} & \frac{-1}{\tilde C} & 0 & 0 & 0 & 0\\ 
\frac{1}{L} & \frac{-R \tilde R}{L \left(R + \mathrm{\bar
      R}\right)} & \frac{-\tilde R}{L \left(R + \mathrm{\bar
      R}\right)} & 0 & 0 & 0\\
 0 & \frac{\tilde R}{\tilde C \left(R + \tilde R\right)} &
 \frac{-1}{\tilde C \left(R + \tilde R\right)} & \frac{-1}{\bar
   C} & 0 & 0\\
 0 & 0 & \frac{1}{L} & \frac{-R \tilde R}{L \left(R +
     \tilde R\right)} & \frac{-\tilde R}{L \left(R +
     \tilde R\right)} & 0\\ 
0 & 0 & 0 & \frac{\tilde R}{\tilde C \left(R + \mathrm{\bar
      R}\right)} & \frac{-1}{\tilde C \left(R + \tilde R\right)}
& \frac{-1}{\tilde C}\\ 
0 & 0 & 0 & 0 & \frac{1}{L} & \frac{-\tilde R}{L} 
\end{array}\right]\\
N&=\left[\begin{array}{cccccc} 0 & 0 & 0 & 0 & 0 & 0\\ 1 & \frac{-R
      \tilde R}{R + \tilde R} & \frac{-\tilde R}{R + \tilde R} & 0 & 0 & 0\\ 0
    & 0 & 0 & 0 & 0 & 0\\ 0 & 0 & 1 & \frac{-R \tilde R}{R + \tilde R} &
    \frac{-\tilde R}{R + \tilde R} & 0\\ 0 & 0 & 0 & 0 & 0 & 0\\ 0 & 0 & 0
    & 0 & 1 & - \tilde R \end{array}\right]\\
B&=\left[\begin{array}{ccccccc} \frac{1}{\tilde C R} & 0 & 0 & 0 & 0 &
    0 & 0 \end{array}\right]^T\\
C&=\left[\begin{array}{ccccccc} -\frac{1}{R} & 0 & 0 & 0 & 0 & 0 & 0 \end{array}\right].
\end{align*}}
For larger $n$, the band structure of $A$ and $N$ is extended periodically.
To see the behaviour of our two methods, we reduce from order $n=20$
to the orders $r=1,3,5,\ldots,19$, and compute both the theoretical bounds
and the actual approximation errors in the $H^\infty$-norm. The
results are shown in the following figure.
 \begin{center}
  \pgfimage[width=.7\linewidth]{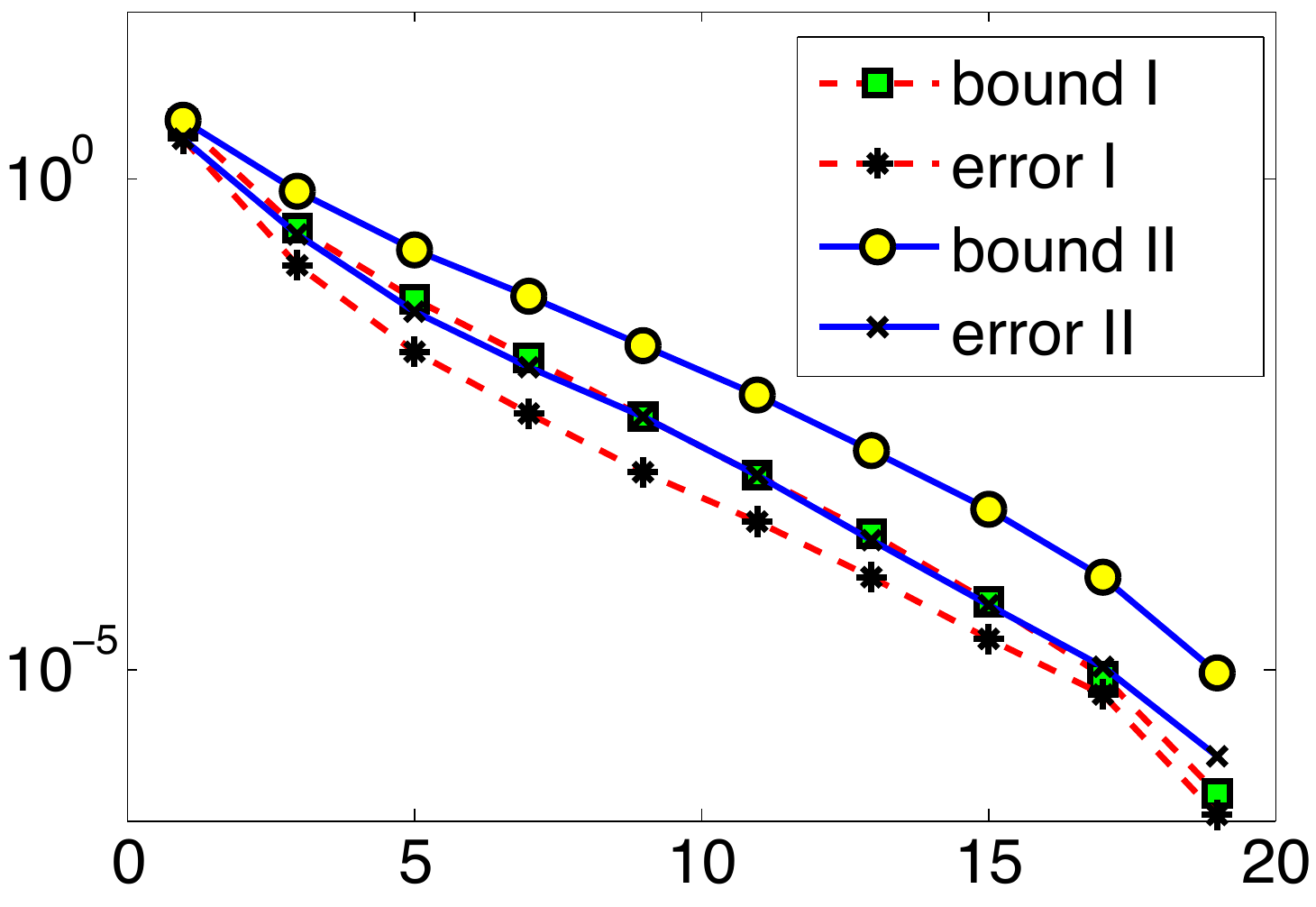}
\end{center}
In this example, for  both types the bounds hold, and for all
reduced orders, type I gives a better approximation than type II.

\subsection{A heat transfer problem}\label{sec:heat}
As another example we consider a stochastic modification of the heat transfer problem described in
\cite{BennDamm11}. On the unit square $\Omega=[0,1]^2$ the heat equation $x_t=\Delta x$
is given with Dirichlet  condition $x=u_j$, $j=1,2,3$ on three
of the boundary edges and a stochastic Robin condition $n\cdot \nabla x=(1/2 +
\dot w)x$ on the fourth edge (where $\dot w$ stands for white noise). A standard 5-point finite difference
discretization on a $10\times 10$ grid
leads to a modified Poisson matrix $A\in\mathbb{R}^{100\times 100}$ and corresponding matrices $N\in\mathbb{R}^{100\times 100}$
and $B\in\mathbb{R}^{100\times 3}$. 
We use the input $u\equiv\left[
  \begin{smallmatrix}
    1\\1\\1
  \end{smallmatrix}
\right]$ and choose the average
temperature as the output, i.e.\
$C=\tfrac{1}{100}[1,\ldots,1]$. We apply balanced
truncation of type I and type II. 
For type II, an LMI-solver (MATLAB\textsuperscript{\textregistered}  function
\texttt{mincx}) is used to compute $P$ as a solution of
the LMI \eqref{eq:lmi} which minimizes $\tr P$ or
$\tr PQ$.

 In the following two figures, we compare the  reduced systems of order
 $r=20$ for
both types. 
The left figure shows the decay of the singular values. 
Since the LMI-solver was called with tolerance level $10^{-9}$,
 only the first about $25$ singular values for type II have the
 correct order of magnitude.
The right figure shows the approximation error
$\|y(t)-y_r(t)\|$ over a given time interval. For both types it has
the same order of magnitude.
 In fact, for many examples we have observed both methods
to yield
very similar results. 

\begin{center}
  \pgfimage[width=.49\linewidth]{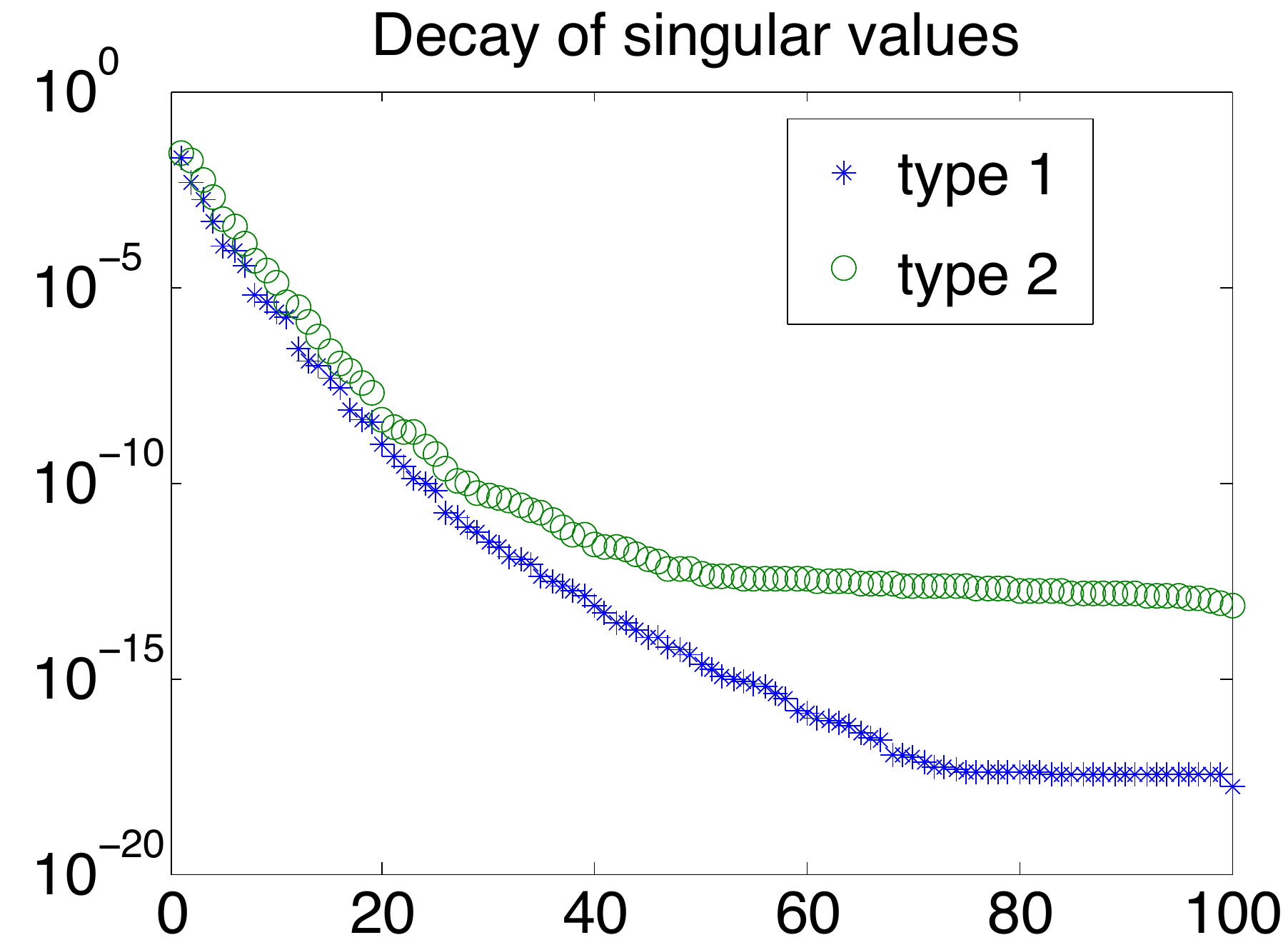}
\pgfimage[width=.49\linewidth]{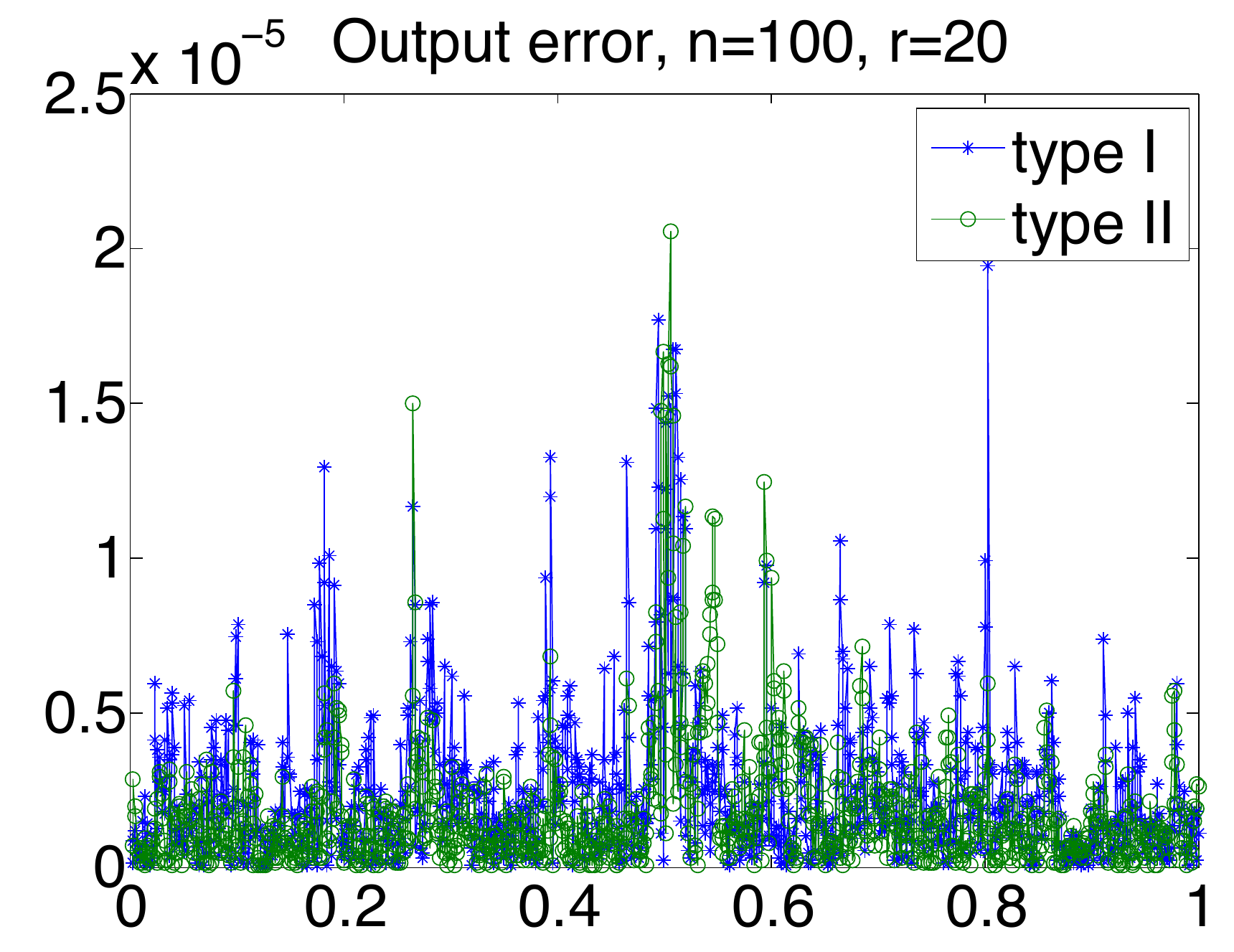}
\end{center}


We have computed the estimated error norm 
and the actual
approximation error for both types: 
\begin{center}
  \begin{tabular}{c||c|c||c|c}
    &$\sum_{j=11}^{100}\sigma_j$&$\|\mathbb{L}-\mathbb{L}_{10}\|$&$\sum_{j=21}^{100}\sigma_j$&$\|\mathbb{L}-\mathbb{L}_{20}\|$\\\hline
    I &$4.66e-06$&$9.30e-06$&$2.00e-09$&$9.65e-09$\\\hline
    II &$1.75e-05$&$4.83e-06$&$1.72e-08$&$9.70e-09$
  \end{tabular}
\end{center}
As we can see, the upper error bound fails for type I, but is correct for type
II. Nevertheless, judging from the $H^\infty$ error, neither of the
types seems to be preferable over the other.

\subsection{Summary}
\label{sec:summary}

Clearly, higher dimensional examples are required to get more insight.
To this end a more sophisticated method for the solution of
\eqref{eq:lmi} is needed. With  general purpose 
LMI-software on a standard
Laptop, we hardly got higher than $n=100$.

\section{Conclusions}

We have compared two types of balanced truncation for stochastic
linear systems, which are related to different Gramian type matrices
$P$ and $Q$. The following table collects properties of these
reduction methods. 
\begin{center}
  \begin{tabular}{l||c|c}
    Type&I&II\\\hline\hline
Def.\ of $P,Q$&\eqref{eq:Gram1}&\eqref{eq:Gram2}\\\hline 
Stability?&Yes, \cite{BennDamm14c}&Yes, Thm.\ \ref{satz:Trunc_Stab2}\\\hline
$H^2$-bound?& Yes, \cite{RedmBenn14}&no result\\\hline
 $H^\infty$-bound?&No, Ex.\ \ref{ex:noerrbound} &Yes, Thm.\
    \ref{thm:NewGramianBound} or \cite{DammBenn14}
  \end{tabular}
\end{center}
The main contributions of this paper are the preservation of
asymptotic stability for type II balanced truncation proved in   Theorem \ref{satz:Trunc_Stab2}
and the new proof of the $H^\infty$ error bound in Theorem \ref{thm:NewGramianBound}. The efficient
solution of \eqref{eq:Gram2.P} is an open issue and requires further
research. The same is true for the computation of the stochastic $H^\infty$-norm.






\appendix

\section*{Asymptotic mean square stability}

Consider the stochastic linear system of It\^o-type
\begin{align}\label{eq:stochNj}
  dx&=Ax\,dt+Nx\,dw\;,
\end{align}
where $w=(w(t))_{t\in\mathbb{R}_+}$ is a zero mean real 
Wiener process on a 
probability space $(\Omega ,{\cal F} ,\mu )$ with respect to an
increasing  family 
$({\cal F}_t)_{t\in \mathbb{R}_+}$ of $\sigma$-algebras ${\cal F}_t
\subset{\cal F} $ (e.g.~\cite{Arno74, Oeks98}).\\
Let $L^2_w(\mathbb{R}_+,\mathbb{R}^q)$ denote the corresponding space of
non-anti\-cipating  stochastic processes $v$ with values in $\mathbb{R}^q$
and norm
\[
\|v(\cdot)\|^2_{L^2_w}:={\cal E}\left(\int_0^\infty\|v(t)\|^2dt\right)<\infty,
\]
where ${\cal E}$ denotes expectation. By definition, system \eqref{eq:stochNj}
is asymptotically mean-square-stable,
if $\mathcal{E}(\|x(t)\|^2)\stackrel{t\to\infty}\longrightarrow0$,
for all initial conditions $x(0)=x_0$.

We have the following version of  Lyapunov's matrix theorem, see
\cite{Khas80}. Here $\otimes$ denotes the Kronecker product.

\begin{satz}\label{thm:Hans} 
The following are equivalent.
  \begin{description}
\item[(i)] System \eqref{eq:stochNj} is asymptotically mean-square
  stable. 
\item[(ii)] $\max\{\Re\lambda\;\big|\; \lambda\in\sigma(A\otimes I+I\otimes A+N\otimes N)\}<0$
\item[(iii)] $\exists Y>0:\exists X>0$: $A^TX+XA+N^TXN=-Y$
\item[(iv)] $\forall Y>0:\exists X>0$: $A^TX+XA+N^TXN=-Y$
\item[(v)] $\forall Y\ge 0:\exists X\ge 0$: $A^TX+XA+N^TXN=-Y$
  \end{description}
\end{satz}

\begin{remark}\rm
  The theorem (like all other results in this paper) carries over to systems  
\begin{align*}
  dx&=Ax\,dt+\sum_{j=1}^kN_jx\,dw_j
\end{align*}
with more than one noise term, and many more equivalent criteria can be
provided, see e.g.\ \cite{Schn65} or \cite[Theorem 3.6.1]{Damm04}.\\
\end{remark}

The following
theorem does not require any stability assumptions (see
\cite[Theorem 3.2.3]{Damm04}). It is central in the analysis of mean-square stability.
\begin{satz}\label{thm:KreinRutman}
Let $$\alpha=\max\{\Re\lambda\;\big|\; \lambda\in\sigma(A\otimes I+I\otimes A+N\otimes N)\}\;.$$
Then there exists a nonnegative
definite matrix $V\neq 0$, such that $$(\mathcal{L}_A^*+\Pi_N^*)(V)=AV+VA^T+NVN^T=\alpha V\;.$$
\end{satz}
We also note a simple consequence of this theorem \cite[Corollary 3.2]{BennDamm14c}. Here
$\langle Y,V\rangle=\tr (YV)$ is the Frobenius inner product for
symmetric matrices.
\begin{cor}\label{cor:semidefY}
Let $\alpha, V$ as in the theorem. 
 For given $Y\ge 0$ assume that
  \begin{align}\label{eq:semidefY}
\exists X>0:\;    \mathcal{L}_A(X)+\Pi_N(X)&\le-Y\;.
  \end{align}
Then $\alpha\le 0$. Moreover, if $\alpha=0$ then $YV=VY=0$.
\end{cor}
\section*{The stochastic bounded
  real lemma}
Now let us consider system \eqref{eq:ls} with input $u$ and output
$y$. If system \eqref{eq:stochNj} is asymptotically mean-square
stable, then \eqref{eq:ls} defines an input output operator
$\mathbb{L}:u\mapsto y$ from $L^2_w(\mathbb{R},\mathbb{R}^m)$ to
$L^2_w(\mathbb{R},\mathbb{R}^p)$, see \cite{HinrPrit98}. By
$\|\mathbb{L}\|$ we denote the induced operator norm, which is an
analogue of the deterministic $H^\infty$-norm. It can be characterized
by the stochastic bounded real lemma. 
\begin{satz}\cite{HinrPrit98}\label{thm:sbrl}
  For $\gamma>0$, the following are equivalent.
  \begin{itemize}
\item[(i)] System \eqref{eq:stochNj} is asymptotically mean-square stable and $\|\mathbb{L}\|<\gamma$.
  \item[(ii)] There exists a negative definite solution $X<0$ to the Riccati
    inequality
\begin{align*}
A^TX+XA+N^TXN-C^TC-\gamma^{-2}XBB^TX>0\;.
\end{align*}
\item[(iii)] There exists a positive definite solution $X>0$ to the Riccati
    inequality
\begin{align*}
A^TX+XA+N^TXN+C^TC+\gamma^{-2}XBB^TX<0\;.
\end{align*}
  \end{itemize}
\end{satz}
We have stated the obviously equivalent formulations (ii) and (iii) to
avoid confusion arising from different formulations in the literature.
Under additional assumptions also non-strict versions can be
formulated. The following sufficient criterion is given in
\cite[Corollary 2.2.3]{Damm04} (where also the signs are changed).
Unlike in the previous theorem, here asymptotic mean-square stability
is assumed at the outset.
\begin{satz}\label{thm:sbrl0}
  Assume that  \eqref{eq:stochNj} is asymptotically 
  stable in mean-square. If there exists a nonnegative definite matrix $X\ge 0$,
  satisfying 
\begin{align*}
A^TX+XA+N^TXN+C^TC+\gamma^{-2}XBB^TX\le 0\;,
\end{align*}
then $\|\mathbb{L}\|\le \gamma$.
\end{satz}


\end{document}